\documentclass[12pt,leqno]{amsart}
\usepackage{amsfonts}                   
\usepackage{latexsym}                   
\usepackage{amssymb}

\newtheorem{teorema}{Theorem}[section] 
\newtheorem{lemma}[teorema]{Lemma}
\newtheorem{propos}[teorema]{Proposition}
\newtheorem{corol}[teorema]{Corollary}
\newtheorem{ex}{Example}[section]
\newtheorem{rem}{Remark}[section]
\newtheorem{defin}[teorema]{Definition}

\def\defin{\par\ifdim\lastskip<\smallskipamount\removelastskip
  \smallskip\fi\noindent{\bf\ignorespaces Definition\unskip:\enspace}\rm
\ignorespaces}

\def\bit{\begin{itemize}}
\def\eit{\end{itemize}}
\def\be{\begin{equation}}
\def\ee{\end{equation}}
\def\beq{\begin{eqnarray}}
\def\eeq{\end{eqnarray}}

\def\ba{\begin{array}}
\def\ea{\end{array}}
\def\bt{\begin{teorema}}
\def\et{\end{teorema}}
\def\bp{\begin{propos}}
\def\ep{\end{propos}}
\def\bl{\begin{lemma}}
\def\el{\end{lemma}}
\def\bc{\begin{corol}}
\def\ec{\end{corol}}
\def\br{\begin{rem}\rm}
\def\er{\end{rem}}
\def\bex{\begin{ex}\rm}
\def\eex{\end{ex}}
\def\bd{\begin{defin}}
\def\ed{\end{defin}}
\def\demo{\par\noindent{\bf Proof.\ }}
\def\enddemo{\ $\Box$\par\vskip.6truecm}

     \def\nin{\noindent}

\def\ES{\varnothing}    

\def\R{{\mathbb {R}}}   \def\a {\alpha} \def\b {\beta}\def\g{\gamma}
\def\N{{\mathbb N}}
\def\B{{\mathbb B}}       \def\e{\varepsilon}
\def\C{{\mathbb C}}
\def\D{{\mathbb D}}      
					\def\p{\partial}
					\def\t{\theta}\def\z{\zeta}

\def\O{\Omega}					
\def\ES{\varnothing}    

  \def\smi{\smallsetminus}

\def\sbs{\!\subset\!}

\def\til{\tilde}
\def\Til {\widetilde}
\def \Hat {\widehat}

\def\oli{\overline}				 \def\tms{\times}

\def\benu{\begin{enumerate}} \def\eenu{\end{enumerate}}
\def\beqn{\begin{eqnarray*}}  \def\eeqn{\end{eqnarray*}}

\def\beqn{\begin{eqnarray*}}  \def\eeqn{\end{eqnarray*}}

\title[Boundary problem for Levi flat graphs]{Boundary problem for Levi flat graphs}
\author[P. Dolbeault, G. Tomassini* and D. Zaitsev**]{Pierre Dolbeault, Giuseppe Tomassini* and Dmitri Zaitsev**}
\address{Giuseppe Tomassini: Scuola Normale Superiore, Piazza dei Cavalieri 7, 56126 Pisa, ITALY}
\email{g.tomassini@sns.it}
\address{Pierre Dolbeault: Institut de Math\'ematiques de Jussieu, Universit\'e
Pierre et Marie Curie, 175 rue du Chevaleret, 75013 Paris, FRANCE}
\email{pierre.dolbeault@upmc.fr}
\address{Dmitri Zaitsev: School of Mathematics, Trinity College Dublin, Dublin 2, Ireland}
\email{zaitsev@maths.tcd.ie}
\thanks{$^{*}$Supported by the project MIUR "Geometric Properties of Real and
Complex Manifolds".}
\thanks{$^{**}$Supported in part by the Science Foundation Ireland grant 06/RFP/MAT018.}
\begin{document}
\numberwithin{equation}{section}
\date{\today}

\begin{abstract} 
In \cite{DTZ2} the authors provided general conditions
on a real codimension $2$ submanifold $S\subset\C^{n}$, $n\ge 3$,
such that there exists a possibly singular Levi-flat hypersurface $M$ bounded by $S$.

In this paper we consider the case when $S$ is a graph of a smooth function over 
the boundary of a bounded strongly convex domain $\Omega\subset\C^{n-1}\times\R$
and show that in this case $M$ is necessarily a graph of a smooth function over $\Omega$.
In particular, $M$ is non-singular.
\end{abstract}

\maketitle

\section{Introduction}

The problem of finding a Levi-flat hypersurface $M\subset\C^n$  with
prescribed boundary $S$ (the complex analogue of the real Plateau's problem), 
has been extensively studied for $n=2$ 
(cf. \cite{Bi, BeG, BeK,  Kr, ChS, Shc, ST, ShT}).
In \cite{DTZ2} (announced in \cite{DTZ1})
we addressed this problem for $n\ge 3$, where the situation is substantially different.
In contrast to the case $n=2$, for $n\ge 3$ the boundary $S$
has to satisfy certain compatibility conditions.
Assuming those necessary conditions as well as 
the existence of complex points, their ellipticity and non-existence
of complex subvarieties in $S$, we have constructed in \cite{DTZ2}
a (unique but possibly singular) solution to the above problem.
An example was also provided in \cite{DTZ2}
showing that one may not always expect a smooth solution $M$ in general.

The purpose of this paper is to show that the solution $M$ is smooth
if the given boundary has certain ``graph form''.
More precisely, in the coordinates $(z,u+iv)\in\C^{n-1}\times\C$,
we assume that $S$ is the graph of a smooth function $g\colon{\rm b}\O\to\R_v$, where ${\rm b}\O$ is the smooth boundary of a strongly convex bounded domain $\O$ in $\C_z^{n-1}\tms\R_u$ and $S$ satisfies the assumptions of \cite{DTZ2} mentioned above. Let $M$ be the solution given by these theorems.
Recall that it is obtained as a projection to $\C^{n}$ of a Levi-flat subvariety with negligible singularities in $[0,1]\times\C^{n}$.
Let $q_1,q_2\in{\rm b}\O$ be the projections of the complex points $p_1,p_2$ of $S$. Using a theorem of Shcherbina on the polynomial envelope of a graph in $\C^2$ (cf. \cite{Shc}) we 
here prove (cf. Theorem \ref{REG}) that
\bit
\item[i)] the solution $M$ is the graph of a Lipschitz function $f\colon\overline\O\rightarrow \R_v$ with $f|_{{\rm b}\O}=g$ which is smooth on $\oli\O\setminus \{q_1,q_2\}$;
\item[ii)] $M_0={\rm graph}(f)\smi S$ is a Levi flat hypersurface in $\C^n$. 
\eit
The regularity of $f$ at $q_{1}$ and $q_{2}$ remains an interesting open problem
closely related to the work of Kenig and Webster \cite{KW1,KW2}.

\section{Preliminaries}\label{PRE}
In this section we collect some facts that will be used in the sequel. 
\subsection{Remarks about Harvey-Lawson theorem}
Let $D$ be a strongly pseudoconvex bounded domain in $\C^n$, $n\ge3$, with boundary ${\rm b}D$, $\Sigma\sbs{\rm b}D$ a compact connected maximally complex $(2d-1)$-submanifold with $d>1$. Then, in view of the theorem of Harvey and Lawson in~\cite[Theorem 10.4]{HL1} (see also~\cite{HL2}), $\Sigma$ is the boundary of a uniquely determined relatively compact subset $V\subset\oli D$ such that: $\oli V\setminus\Sigma$ is a complex analytic subset of $D$ with finitely many singularities of pure dimension $d$ and, near $\Sigma$, $\oli V$ is a $d$-dimensional complex manifold with boundary. We refer to $V=V_\Sigma$ as the {\it solution of the boundary problem corresponding to $\Sigma$}. A simple consequence is the following:
\bl\label{INTER}
Let $D\subset\C^{n}$ be as above and $\Sigma_1$, $\Sigma_2$ connected, maximally complex $(2d-1)$-submanifolds of ${\rm b}D$. Let $V_1$, $V_2$ be the corresponding solutions of the boundary problem. If $d>1$, $2d>n$ and $\Sigma_1\cap\Sigma_2=\ES$, then $V_1\cap V_2=\ES$.   
\el
\demo
Suppose $V_1\cap V_2\ne\ES$.
Then $2d>n$ implies $\dim V_{1}\cap V_{2}\ge1$.
Since $V_{1}\cap V_{2}$ is an analytic subset of $D$, its closure
 $\oli{V_{1}\cap V_{2}}$ must intersect ${\rm b} D$ and hence also 
$\Sigma_1\cap\Sigma_2\ne\ES$, which contradicts the assumption.
\enddemo 
\subsection{Known results}\label{known}
First, we have the following: a real $2$-codimensional submanifold $S$ of $\C^{n}$, $n\ge 3$, which locally bounds a Levi flat hypersurface must be nowhere minimal near a {\rm CR} point, i.e.\ all local CR orbits must be of positive codimension (cf. \cite[Section 2]{DTZ2}). If $p\in S$ is a complex point, consider local holomorphic coordinates  
$(z,w)\in \C_z^{n-1}\times \C_w$, vanishing at $p$, such that $S$ is locally given by the equation
\begin{equation}\label{FORM} 
w= Q(z)+O(|z|^3),
\end{equation}
where $ Q(z)$ is a complex valued quadratic form in the real coordinates $({\sf Re}\,z, {\sf Im}\,z)\in\R^{n-1}\times\R^{n-1}$. Observing that not all quadratic forms $Q$ can appear when $S$ bounds a Levi flat hypersurface one comes to the condition that $p$ must be {\it flat},
i.e.\ $Q(z)\in {\mathbb R}$ in suitable coordinates.
A natural stronger condition is that of {\it ellipticity}
which means by definition that $Q(z)\in {\mathbb R_+}$ for every $z\ne0$
in suitable coordinates.

Assume that:
\begin{enumerate}
\item $S$ is compact, connected and nowhere minimal at its {\rm CR} points;
\item $S$ has at least one complex point and every such point of  is flat and elliptic;
\item $S$ does not contain complex manifold of dimension $(n-2)$.
\end{enumerate}
Then in \cite[Proposition 3.1]{DTZ2} it was proved  that
\bit 
\item[a)] $S$ is diffeomorphic to the unit sphere with two complex points $p_1, p_2$;
\item[b)] the CR orbits of $S$ are topological $(2n-3)$-spheres that can be represented as level sets of a smooth function $\nu:S\to\R$, inducing on $S_0=S\setminus\{p_1, p_2\} $ a foliation $\mathcal F$ of class $C^\infty$ with $1$-codimensional compact leaves. 
\eit

Next,  by applying a parameter version of Harvey-Lawson's theorem \cite[Theorem~8.1]{HL1}, we obtained in \cite[Theorem 1.3]{DTZ2} a solution to the boundary problem as follows:
\bt\label{main}
Let $S\subset \C^n$, $n\ge 3$ satisfy the above conditions. Then there exist a smooth submanifold $\Til S$ and a Levi flat $(2n-1)$-subvariety $\widetilde M$ in $\C^n\times[0,1]$ (i.e. $\Til M$ is Levi flat in $\C^n\times\C$) such that $\Til S=d\Til M$ in the sense of currents and the natural projection $\pi\colon \C^n\tms[0,1]\to \C^n$ restricts to a diffeomorphism between $\Til S$ and $S$.
\et
As for the singularities of $\Til M$ we have the following results \cite[Theorems 1.4]{DTZ2}:
\bt\label{precise}
The Levi-flat $(2n-1)$-subvariety $\Til M$ can be
chosen with the following properties:
\begin{enumerate}
\item[(i)] $\Til S$ has two complex points $\til p_0$ and $\til p_1$
with $\Til S\cap(\C^n\times\{j\})=\{\til p_j\}$ for $j=0,1$;  every other slice $\C^n\tms\{x\}$ with $x\in(0,1)$, intersects $\Til S$ transversally along a submanifold diffeomorphic to a sphere that bounds (in the sense of currents) the (possibly singular) irreducible complex-analytic hypersurface $(\Til M\setminus\Til S)\cap (\C^n\times\{x\})$;
\item[(ii)] the singular set ${\rm Sing}\,\Til M$ is the union of $\Til S$ and a 
closed subset of $\Til M\setminus\Til S$ 
of Hausdorff dimension at most $2n-3$; moreover each slice
$({\rm Sing}\,\Til M\setminus\Til S)\cap (\C^n\times\{x\})$
is of Hausdorff dimension at most $2n-4$;
\item[(iii)] there exists a closed subset $\Til A\subset\Til S$ 
of Hausdorff $(2n-2)$-dimensional measure zero
such that away from $\Til A$, $\Til M$ is a smooth submanifold with 
boundary $\Til S$ near $\Til S$;
moreover $\Til A$ can be chosen such that each slice $\Til A\cap (\C^n\times\{x\})$
is of Hausdorff $(2n-3)$-dimensional measure zero.
\end{enumerate}
\et
\section{The case of graph}
From now on we assume that $S\subset\C^{n}$, $n\ge 3$, is a graph.
Consider $\C^{n}=\C^{n-1}_z\tms\C_w$ with complex coordinates  $z=(z_1,\ldots,z_{n-1})$ and $w$ where $z_\a=x_\a + iy_\a$, $1\le\a\le n-1$, $w=u+iv$. We also write $\C^n=\big(\C^{n-1}_z\tms\R_u\big)\tms i\R_v$. Accordingly, a point of $\C^{n}$ will be denoted by $(z, u,v)=(z,u+iv)$.

Let $\O$ be a bounded strongly convex domain of $\C^{n-1}_{z}\times\R_u$
with smooth boundary ${\rm b}\O$. By strong convexity here we mean that the second fundamental form of the boundary ${\rm b}\O$ of $\O$
is everywhere positive definite. In particular, $\O\tms i\R_v$ is a strongly pseudoconvex domain in $\C^n$.

Let $g:{\rm b}\O\rightarrow \R_{v}$ be a smooth function, and $S\subset \C^n$ the graph of $g$. We assume that $S$ satisfies the conditions of \cite[Theorem 1.3]{DTZ2} and denote $q_1,q_2\in{\rm b}\O$ the natural projections of the complex points $p_1,p_2$ of $S$, respectively.

Our goal is to prove the following:
\bt\label{REG} 
Let $q_1,q_2\in{\rm b}\O$ be the projections of the complex points $p_1,p_2$ of $S$, respectively. Then, there exists a Lipschitz function $f:\oli\O\rightarrow \R_v$ which is smooth on $\oli\O\setminus \{q_1,q_2\}$ and such that $f_{\vert {\rm b}\O}=g$ and $M_0={\rm graph}(f)\smi S$ is a Levi flat hypersurface of $\C^n$. Moreover, each complex leaf of $M_0$ is the graph of a holomorphic function $\phi\colon\O'\to \C$ where $\O'\subset\C^{n-1}$
is a domain with smooth boundary (that depends on the leaf) and $\phi$ is smooth on $\oli\O'$.
\et
The natural candidate to be the graph $M$ of $f$ is $\pi\big(\Til M\big)$ where $\Til M$ and $\pi$ are as in Theorem \ref{main}. We prove that this is the case proceeding in several steps. 
\subsection{Behaviour near $S$}
Set $m_1=\min\limits_S g$, $m_2=\max\limits_S g$ and $r\gg 0$ such that 
$$
D=\O\times[m_1,m_2]\Subset\B(r)\cap\big(\O\times i\R_{v}\big)
$$
where $\B(r)$ is the ball $\{\vert(z,w)\vert<r\}$. 

Let $\Sigma$ be a CR-orbit of the foliation of $S\setminus\{p_1,p_2\}$. Then, $\Sigma$ is a compact maximally complex $(2n-3)$-dimensional
real submanifold of $\C^n$, which is contained in the boundary of the strongly pseudoconvex domain  $\O\times i\R_v$ of $\C^n$.  Let $V$ be the solution to the boundary problem corresponding to $\Sigma$, i.e.\ the complex-analytic subvariety of  $\O\times i\R_v$ bounded by $\Sigma$. We refer to $V$ as the {\em leaf} bounded by $\Sigma$. From Theorems~\ref{main} and \ref{precise} it follows that $V$ is obtained as projection $\pi(\Til V)$, where $\Til V= (\Til M\setminus\Til S)\cap (\C^{n}\times\{x\})$
for suitable $x\in (0,1)$. In particular, if $M$ denotes $\pi(\Til M)$, $\pi_{|\Til V}$ defines 
a biholomorphism $\widetilde {V}\simeq V$ and 
$M\setminus S\subset D$.

Now let $\Sigma_1$ and $\Sigma_2$ be two distinct {\rm CR} orbits of 
the foliation of $S\setminus\{p_1,p_2\}$, and let $\oli V_1$, $\oli V_2$ be the corresponding leaves
bounded by them. Then $\oli V_1$, $\oli V_2$ do not intersect by Lemma~\ref{INTER}.
\br
In the previous discussion, we only employed the fact that
$\O\times \R_v$ is a strongly pseudoconvex domain and $S$ is
contained in its boundary, without regarding the graph nature of
$S$.  It can happen that the leaves have isolated singularities.
We shall show that this cannot happen in our case.
\er
\bl\label{local}
Let $p\in S$ be a CR point.
Then, near $p$, $M$ is the graph of a function $\phi$ on a domain $U\sbs\C_z^{n-1}\tms\R_u$, which is smooth up to the boundary of $U$.
\el
\demo
Near $p$, $S$ is foliated by local CR orbits. As a consequence of Theorem~\ref{main}, each local CR orbit extends to a compact global CR orbit $\Sigma$ that bounds a complex codimension $1$ subvariety $V_{\Sigma}\subset\O\times i\R_{v}$. Furthermore, near $p$, each $\Sigma$ is smooth and can be represented as the graph of a CR function over a strongly pseudoconvex hypersurface and $V_{\Sigma}$ as the graph of the local holomorphic extension of this function. 
It follows from the Hopf Lemma that $V$ is transversal to the strongly pseudoconvex hypersurface $b\O\times i\R_{v}$ near $p$. Hence the family of $V_{\Sigma}$ near $p$ forms a smooth real hypersurface with boundary on $S$ that can be seen as the graph of a smooth function $\phi$ from a relative open neighbourhood $U$ of $p$ in $\overline\O$ into $\R_{v}$. Finally, Lemma~\ref{INTER} guarantees that this family does not intersect any other leaf $V$ from $M$. This completes the proof.
\enddemo

\bc\label{loc-leaf}
If $p\in S$ is a CR point, each complex leaf $V$ of $M$, near $p$, 
is the graph of a holomorphic function on a domain $\O_{V}\subset\C^{n-1}_{z}$,
which is smooth up to the boundary of $\O_{V}$.
\ec
\demo
Since $M$ is the graph of a smooth function near $p$, its tangent space at every point near $p$ is transversal to 
$i\R_{v}$. Hence the complex tangent space of $M$ at every point near $p$ is transversal to $\C_{w}$.
Since the tangent spaces of the complex leaves of $M$ coincide with the complex
tangent spaces of $M$, it follows that each leaf $V$ projects immersively to $\C^{n-1}_{z}$
and the conclusion follows.
\enddemo

\subsection{$M$ is the graph of a Lipschitz function.}
Assume as before that $\O$ is strongly convex. We have the following
\bp\label{graph}
$M$ is the graph of a Lipschitz function $f:\oli\O\to\R_v$.
\ep
\demo
We fix a nonzero vector $a\in \C^{n-1}_{z}$ and for a given point $(\z,\xi)\in\O$ denote by
$H_{(\z,\xi)}\subset\C_z^{n-1}\tms\{\xi\}$ the complex line through $(\z,\xi)$ in the direction of $(a,0)$.
Furthermore, we set
$$
L_{(\z,\xi)}=H_{(\z,\xi)}+\R(0,1), \quad
\O_{(\z,\xi)}=L_{(\z,\xi)}\cap \O, \quad
S_{(\z,\xi)}=\big(H_{(\z,\xi)}+\C(0,1)\big)\cap S
$$
Then $S_{(\z,\xi)}$ is contained in the strongly convex cylinder
$$
\big(H_{(\z,\xi)}+\C(0,1)\big)\cap \big({\rm b}\O\tms i\R_v\big)
$$
over $H_{(\z,\xi)}+\C(0,1)\simeq\C^2$ and it is the graph of $g_{\vert{\rm b}\O_{(\z,\xi)}}$. 

Since $\O_{(\z,\xi)}=\O\cap L_{(\z,\xi)}$, in view of the main theorem of \cite{Shc}, the polynomial hull $\Hat S_{(\z,\xi)}$ of $ S_{(\z,\xi)}$ is a continuous graph over ${\overline\O}_{(\z,\xi)}$. Consider $M=\pi(\widetilde M)$ and set 
$$
M_{(\z,\xi)}= (H_{(\z,\xi)}+\C(0,1))\cap M.
$$ 
Since $M$ is a union of irreducible analytic subvarieties of codimension $1$ in $\C^{n}$ with boundary in the graph $S$,
each intersection $M_{(\z,\xi)}$ is the union of a
family $\mathcal A$ of $1$-dimensional analytic subsets. 
Clearly, the boundary of a connected component of any such analytic set is contained in
$S_{(\z,\xi)}$. It follows that $M_{(\z,\xi)}$ is contained in the
polynomial hull $\Hat S_{(\z',\xi)}$ of $ S_{(\z,\xi)}$. In view of the main theorem of Shcherbina \cite{Shc}, $\Hat S_{(\z,\xi)}$ is a graph over ${\overline\O}_{(\z,\xi)}=\overline\O\cap L_{(\z,\xi)}$, foliated by analytic discs, so $M_{(\z,\xi)}$ is a graph over a subset $U$ of ${\overline\O}_{(\z,\xi)}$. 

On the other hand, every analytic disc $\Delta$ of $\Hat S_{(\z,\xi)}$ 
has its boundary on $S_{(\z,\xi)}\subset S$.
Since all elliptic complex points are isolated, the boundary of $\Delta$
contains a CR point $p$ of $S$. In view of Lemma~\ref{local}, near $p$,
$M_{(\z,\xi)}$ is also a graph over ${\overline\O}_{(\z,\xi)}$.
Thus, near $p$, we must have $M_{(\z,\xi)}=\Hat S_{(\z,\xi)}$.
In particular, near $p$, $\Delta$ is contained in $M_{(\z,\xi)}$,
and therefore in a leaf $V_\Sigma$ of $M$.
Since $V_\Sigma$ is a closed analytic subset in $\C^{n}\setminus S$,
the whole disc $\Delta$ is contained in $V_\Sigma$ and hence in $M$.
Moreover, $\Delta\subset H_{(\z,\xi)}+\C(0,1)$ thus
we conclude that $\Delta\subset M_{(\z,\xi)}$.
Therefore, every analytic disc of $\Hat S_{(\z,\xi)}$ is contained in $M_{(\z,\xi)}$,
consequently $M_{(\z,\xi)}$ and $\Hat S_{(\z,\xi)}$ coincide. 
It follows that $M$ is the graph of a function $f:\oli\O\to\R_u$. 

Let us prove that $f$ is a continuous function. Choose $(\z,\xi)\in \O$ and a complex line $H_{(\z,\xi)}$ as before. Consider a neighborhood $U$ of $(\z,\xi)$ in $\C_z^{n-1}\tms\R_u$. For $q\in U$, let $H_q$ be the translated of $H_{(\z,\xi)}$ which passes through $q$. With the notation corresponding to the one employed above, we can state the following. For a small
enough neighborhood $V\subset U$ of $p$ in $\C_z^{n-1}\tms\R_u$, let $\Hat S_q$ be the
polynomial hull of $S_q$ in $H_q + \C(0,1)$, and let
$$
\mathcal S_U= \bigcup_{q\in U} \Hat S_q;
$$
then $\mathcal S_U$ is the graph of a continuous function. Indeed let $\oli q$ be a point in $V$, and let $\{q_m\}_{m\in\N}$ be
a sequence of points such that $q_n\to \oli q$. Then, obviously, the
sets $S_{q_m}$ converge to the set $S_{\oli
q}$ in the Hausdorff metric as $n\to \infty$. Moreover, it is also
clear that $\widetilde \O_{q_n} \to \widetilde \O_{\oli q}$ for
$n\to \infty$. Then, by \cite[Lemma 2.4]{Shc} it follows that
$\Hat S_{q_m}\to \Hat S_{\oli q}$ as $m\to\infty$. Since every
$\Hat S_q$ is a continuous graph, this allows to prove easily
that $\mathcal S_U$ is a continuous graph as a whole.

Thus, $f$ is continuous on $\O$, whence on $\oli\O\smi\{q_1,q_2\}$ in view of Lemma \ref{local}. Continuity at $q_1$ is proved as follows. Let let $\{a_m\}_{m\in\N}\sbs\O$ be a sequence of points which converges to $q_1$. Each point $\big(a_m,f(a_m)\big)$ belongs to a complex leaf $V_{\Sigma_m}$ of 
$M$ which is bounded by a compact {\rm CR} orbit $\Sigma_m$ of the foliation of $S\smi\{p_1,p_2\}$ (cf. Section \ref{PRE}). By the maximum principle, for every $m\in\N$ there exists a point $(b_m, g(b_m))$ in $\Sigma_m$ such that
$$
\vert (q_1,g(q_1))-(a_m, f(a_m))\vert\le\vert (q_1,g(q_1))-(b_m, g(b_m))\vert.
$$
We claim that 
$$
\vert (q_1,g(q_1))-(b_m, g(b_m))\vert\to 0
$$
as $m\to\infty$. If not there exists an open $B=B(q_1,r)\O\tms\R_u$ centered at $q_1$ such that $b_m\notin\oli B$ for all $m$. It follows
that 
$$
\Sigma_m\cap\pi^{-1}(\oli B)=\varnothing
$$
for all $m$ and
$$
V_{\Sigma_m}\cap\pi^{-1}(B)\neq\varnothing
$$ 
for $m\gg 0$. This violates the Kontinuit\"atsatz since $\O\tms i\R_v$ is a domain of holomorphy.

Continuity at $q_2$ is proved in a similar way. 

Thus $f$ is continuous on $\oli\O$ and smooth near ${\rm b}\O\smi\{q_1,q_2\}$.

In order to show that $f$ is Lipschitz we now observe that, as it is easily proved, $f_{|\O}$ is a weak solution of the {\em Levi-Monge-Amp\`ere} operator defined in \cite{ST} with smooth boundary value, so, in view of \cite[Theorems 2.4, 4.4, 4.6]{ST}, it is Lipschitz. This concludes the proof of Proposition \ref{graph}.
\enddemo
\br
$M$ is the envelope of holomorphy of $S$.
\er
\subsection{Regularity.}
In order to prove that $M\smi\{p_1,p_2\}$ is a smooth manifold with boundary we need the following:
\bl\label{RELF}
Let $U$ be a domain in $\C^{n-1}_z\times \R_u$, $n\ge 2$, $f:U\to \R_v$ a continuous function. Let $A\sbs{\rm graph}(f)$ be a germ of complex analytic set of codimension $1$. Then $A$ is a germ of a complex manifold, which is a graph over $\C^{n-1}_{z}$.
\el 
\demo
The idea of the proof (here is slightly modified) is due to Jean-Marie Lion {cfr. \cite{Lion}.

Let us denote by $z_1,\ldots,z_{n-1},w=u+iv,$ the complex coordinates in $\C^{n-1}_z\tms\C_w$. We may suppose that $A$ is a germ at $0$. Let $h\in\mathcal O_{n+1}$ be a non identically vanishing germ of holomorphic
function such that $A=\{h=0\}$. Let $\D_\e$ be the disc $\{z=0\}\cap\{|w|<\e\}$. Then, for $\e<<1$, we have either $A\cap \D_\e=\{0\}$ or $A\cap \D_\e = \D_\e$. The latter is not possible since $\D_\e$ is not contained in any graph over $\C^{n-1}\times \R_u$. It follows that $A\cap \D_\e=\{0\}$, i.e.\ $A$ is $w$-regular. 
Let us denote by $\pi$ the projection $\C^{n}_z\to\C^{n-1}_z$. Then, by the local parametrization theorem for analytic sets there exists $d\in \N$ such that
\begin{itemize}
\item for some neighborhood $U$ of $0$ in $\C^{n-1}_z$, there exists an analytic set $\Delta\subset U$ such that $A_\Delta=A\cap ( (U\setminus \Delta)\times\D_\e)$ is a manifold;
\item $\pi:A_\Delta\to U\setminus \Delta$ is a $d$-sheeted covering.
\end{itemize}
\nin We claim that the covering $\pi:H_\Delta\to U\setminus \Delta$ is trivial.
Otherwise, there would exist a closed loop $\g:[0,1]\to U\setminus \Delta$
whose lift $\widetilde \g$ to $A_{\Delta}$  is not closed. 
We extend $\g$ to $\R$ by periodicity and extend $\widetilde\g$ to $\R$ as lift of $\g$.
Define 
$\a=u \circ \widetilde \g = u \circ \g$, $\b=v \circ \widetilde \g$.
Since
$\a$ is continuous and bounded, there exists $\t\in \R$ such that 
$\a(\t) =\a(\t +1)$. But then $\b(\t)=\b(\t+1)$ since by the assumption, 
$\b(\t) = f(\g(\t),\a(\t))$. Hence $\widetilde \g(\t)
= \widetilde \g(\t+1)$, a contradiction with the assumption that $\widetilde\g$ is not closed.

Since $\pi:A_\Delta\to U\setminus \Delta$ is a trivial covering, we may
define $d$ holomorphic functions $\tau_1,\ldots,\tau_d:U\setminus
\Delta\to \C$ such that $A_\Delta$ is a union of the graphs of the $\tau_j$'s.
By Riemann's extension theorem, the functions $\tau_j$ extend as holomorphic functions $\tau_j\in \mathcal O(U)$. The desired conclusion will follow from the
fact that all the $\tau_j$ coincide. Indeed, suppose, by contradiction,
$\tau_1\neq \tau_2$; then for some disc $\D\subset U$ centered at $0$
we have $\tau_1|_\D \neq \tau_2|_\D$ and then, after shrinking $\D$,
$(\tau_1-\tau_2)|_\D$ vanishes only at $0$. But, by virtue of the hypothesis, $\{{\sf Re}\,(\tau_1-\tau_2)=0\}\subset \{\tau_1 -\tau_2=0\} =
\{0\}$, when restricted to $\D$.
The latter is not possible since $(\tau_1-\tau_2)|_{\D}\ne0$ is
holomorphic and thus an open map (whose image must include a segment of the imaginary axis).
\enddemo

{\bf Proof of Theorem \ref{REG}}. Consider the foliation on $S\smi\{p_1,p_2\}$ given by the level sets of the smooth function $\nu\colon S\to [0,1]$ as in Section \ref{PRE} and set $L_t=\{\nu=t\}$ for $t\in (0,1)$. Let $V_t\subset\overline\O\times i\R_{v}\subset\C^{n}$ be the complex leaf of $M$ bounded by $L_t$ and $\pi\colon\C^{n-1}_{z}\times\C_{w}\to\C_z^{n-1}$ denote the natural projection. We have:
\bit
\item  by Proposition~\ref{graph}, $M$ is the graph of a continuous function over $\O$ and by Lemma~\ref{RELF}, each leaf $V_{t}$ is a complex hypersurface and $\pi|_{V_{t}}$ is a submersion.
\item Since $\O$ is strongly convex, an argument completely analogous to that of \cite[Lemma~3.2]{Shc} shows that 
 $\pi_{|V_t}$ is one-to-one, then,  by Corollary~\ref{loc-leaf}, $\pi$ sends $V_t$ onto a domain $\O_t\subset\C^{n-1}_{z}$ with smooth boundary.

If 
\bit
\item[] $\pi_u:\big(\C^{n-1}_z\tms\R_u\big)\tms i\R_v\to\R_u,$
\item[]$\pi_v:\big(\C^{n-1}_z\tms\R_u\big)\tms i\R_v\to\R_v$
\eit
denote the natural projections then ${\pi_u}_{|L_t}=a_t\circ\pi_{|L_t}$ and ${\pi_v}_{|L_t}=b_t\circ\pi_{|L_t}$, where $a_t$ and $b_t$ are smooth functions in ${\rm b}\O_t$. Furthermore, the boundary ${\rm b}\O_t$ and $a_{t}$, $b_{t}$
depend smoothly on $t$ for $t\in (0,1)$. 
The latter property means that one has a local parametrization of ${\rm b}\O_t$
smoothly depending on $t$ and such that the functions $a_{t}$, $b_{t}$ also depend smoothly on $t$ when composed with this parametrization. It follows that
\item if $(z_t,w_t)\in M$, then $w_t=u_t+iv_t$ is varying in $V_t$, 
so $u_t+iv_t$ is the holomorphic extension to $\O_{t}$ of $a_{t}+ib_{t}$.
In particular, $u_{t}$ and $v_{t}$ are smooth functions in $(z,t)$,
e.g.\ as a consequence of the Martinelli-Bochner formula.
\item
The derivative $\p u_t/\p t$ 
is defined and harmonic in $\O_{t}$ for each $t$,
and has a smooth extension to the boundary ${\rm b}\O_{t}$.
Moreover, it follows from Lemma~\ref{local} and Corollary~\ref{loc-leaf}
that $\p u_t/\p t$ does not vanish on ${\rm b}\O_{t}$.
Since the CR orbits $L_{t}$ are connected in view of Theorem~\ref{main},
the boundary ${\rm b}\O_{t}$ is also connected and hence $\p u_t/\p t$ 
has constant sign on ${\rm b}\O_{t}$.
Then, by the maximum principle, $\p u_t/\p t$ has constant sign in $\O_{t}$
and, in particular, does not vanish.
The latter implies the $M\setminus S$ is the graph of a smooth function
over $\O$, which extends smoothly to $\overline \O\setminus \{q_{1},q_{2}\}$.
 \item It furthermore follows from Proposition~\ref{graph} that $M$ is the graph of a Lipschitz function over $\overline \O$.
 This completes the proof of Theorem~\ref{REG}.
\eit

\end{document}